\documentclass[10pt,reqno]{article}
\usepackage{latexsym}
\usepackage{amssymb}
\usepackage{amsmath}
\usepackage[all]{xy}

\newtheorem{thm}{Theorem}[section]
\newtheorem{defin}{Definition}[section]
\newtheorem{prop}{Proposition}[section]
\newtheorem{lemma}{Lemma}[section]

\newcommand{\Mbar}{\overline {M}}
\newcommand{\Sbar}{\overline {S}}
\newcommand{\Ric}{\text{Ric}}
\newcommand{\Vol}{\text{Vol}}

\newcommand{\jbar}{\bar {j}}

\newcommand{\zbar}{\bar {z}}
\newcommand{\lbar}{\bar {l}}
\newcommand{\pbar}{\bar {p}}

\newcommand{\kahler}{K\"ahler     }
\newcommand{\KE}{K\"ahler-Einstein     }
\newcommand{\qed}{\hfill \fbox{}}

\def\NN {{\mathbb N}}

\def\RR {{\mathbb R}}
\def\CC {{\mathbb C}}

\def\Om{\Omega}

\def\de{\delta}
\def\om{\omega}

\def\ve{\varepsilon}

\title{On complete Ricci-flat metrics on open \kahler manifolds}

\author{Bianca Santoro}
\date{November 11th, 2006}

\begin{document}

\maketitle

\section{Introduction}
\label{s.introduction}

In 1978, Yau \cite{Y1} proved the Calabi Conjecture, by showing
existence and uniqueness of \kahler metrics with prescribed Ricci curvature
on compact complex manifolds.
Here the complex manifolds in question are already supposed to admit
a \kahler metric whose Ricci form satisfies the natural conditions arising
from Chern-Weil theory.
The remarkable paper \cite{Y1} also
establishes several related results of fundamental importance
in the study of complex manifolds.

The results are proved by reducing them to questions in non-linear partial
differential equations of Monge-Amp\`ere type. Once this is done, the questions
are treated by a continuity method involving intrincated {\em a priori} estimates.
These methods, in fact, are interesting in themselves.

Following the work by Yau, Tian and Yau \cite{TY1} settled a non-compact version
of Calabi's Conjecture on quasi-projective manifolds that can be compactified
by adding a smooth, ample divisor.
Again, the desired metric was constructed by translating the geometric problem
into a complex Monge-Amp\`ere equation. Unfortunately, however,
the solution to the  Monge-Amp\`ere equation obtained in \cite{TY1} is only
shown to be bounded at infinity. An interesting question that arises naturally
is then to work out the asymptotic behavior of these completes metrics
near the divisor (regarded as ``infinity''). This question is, in fact, posed
by Tian and Yau in \cite{TY1}. Clearly the knowledge of full asymptotic expansion
of these metrics immediately yields a more
accurate picture of their behavior near infinity.

In a subsequent work (\cite{TY2}),  Tian and Yau
extended their result for the case where the divisor has
multiplicity greater than one, and is allowed to have orbifold singularities.
A generalization that was done independently by Bando and
Kobayashi \cite{K}. Later, Joyce \cite{Joyce} provided the sharp asymptotics
for the decay of the solutions provided in \cite{TY2} thus clarifying, in this
case, the geometry associated to this metric near infinity.

The problem of computing the asymptotic expasion of the solutions
to the Monge-Amp\`ere equations obtained in \cite{TY1} however remained open,
since the case where the divisor has multiplicity exactly~$1$ contains
special difficulties resisting to the method of Joyce \cite{Joyce}.

In the author's thesis (\cite{Santoro1}, \cite{Santoro2}),
this question is treated in detail and the corresponding asymptotics of the solution
to the Monge-Amp\`ere equation in \cite{TY1} is obtained in two steps.
First, it is  considered a sequence of complete \kahler metrics with special approximating
properties on $M$. Whereas these metrics are only ``approximate solutions'', they
have the advantage of being given by {\it explicit formulas}. The second
step consists then of using
these metrics as ambient metrics to solve the Monge-Amp\`ere equation in \cite{TY1}.
Thanks to the approximating properties of the mentioned metrics, we shall be able
to completely describe the decay of the solutions at infinity.

The present text has two purposes. First it partially extends
the result in \cite{Santoro2}, by considering a larger
class of quasi-projective manifolds, which is described in Section~\ref{s.general}.
Basically the definition of this class consists of relaxing the condition that
the divisor ``at infinity'' should be ample. This will allow us to encompass,
for example, manifolds with cylindrical ends as those considered by \cite{Ko},
cf. Section~3.

On the other hand, this text is also intended to provide
a quick introductory reference to the study of
Ricci-flat metrics on open manifolds. For reasons of space, we shall not
report on the analogous problem for K\"ahler-Einstein metrics except for
a brief comment in Section~4. We point out however that the analysis of
the differential operators involved in each case is very different.

It is assumed that the reader is familiar with
some basic facts in Riemannian and Complex Geometry, as well as with the analysis of
elliptic partial differential equations.

\section{Background: Calabi Conjecture on Compact Manifolds}
\label{s.background}

In order to understand the main problem
considered in this work,
we shall start from its original motivation.

\bigskip

Let $M$ be a compact, complex manifold of complex dimension
$n$, and consider $g$,
a hermitian metric defined on $M$.
Note that $g$ is a complex-valued sesquilinear form acting on
$TM \times TM$, and can therefore be written as
$$
g = S - 2\sqrt{-1} \om_g,
$$
where $S$ and $-\om$ are real bilinear forms.

If $(z_1, \dots, z_n)$ are local coordinates around a point $x\in M$,
we can write the metric $g$ as $\sum g_{i\jbar}dz^i \otimes d\zbar^j$.
Then, it is easy to see that in these coordinates
$$
\om_g = \frac{\sqrt{-1}}{2} \sum_{i, j = 1}^n g_{i\jbar}dz^i \wedge d\zbar^j.
$$
The form $\om_g$ is a real $2-$form of type $(1,1)$, and is called
the {\it fundamental form} of the metric $g$.

\begin{defin}
We say that a hermitian metric on a complex manifold is
{\em \kahler} if its associated fundamental form $\om_g$
is closed, {\em ie}, $d \om_g = 0$.
A complex manifold equipped with a \kahler metric is
called a {\em \kahler manifold}.
\end{defin}

The reader will find in the literature a number of equivalent definitions
for a \kahler metric. We will keep this choice for convenience of the
exposition.
Also,  note that, on a \kahler manifold,
the form $\om_g$ is uniquely determined
by the metric $g$, and vice-versa.

Let $R(g) = R_{i\jbar k \lbar}$ be the Riemann curvature tensor of the metric $g$
written in the above mentioned coordinates.
We define the {\it Ricci curvature tensor} of the metric $g$ as being the trace of the
Riemann curvature tensor. Its components in local coordinates can be written
as
\begin{equation}
\label{e.Riccidef}
\Ric_{k\lbar} =  \sum_{i, j = 1}^n g^{i\jbar} R_{i\jbar k \lbar} =
-\frac{\partial^2}{\partial z_k \partial \zbar_l} \log \det (g_{i\jbar}).
\end{equation}
The{\it Ricci form} associated to $g$ can then be defined by setting
$$
\Ric = \sum_{i, j = 1}^n Ric_{i\jbar} dz^i \wedge d\zbar^j.
$$
in local coordinates.

Now, given a metric $g$, we can define a matrix-valued $2$-form $\Om$ by
writing its expression in local coordinates, as follows
\begin{equation}
\label{e.curvform}
\Om_i^j = \sum_{i, p = 1}^n g^{j\pbar}  R_{i\pbar k \lbar} dz^k \wedge d\zbar^l.
\end{equation}
This expression for $\Om$ gives a well-defined matrix of $(1,1)$-forms, to
be called the {\it curvature } of the metric $g$.

Following Chern-Weil Theory, we want to look at  the following expression
$$
\det\left(\text{Id} + \frac{t\sqrt{-1}}{2\pi} \Om \right) = 1 + t\phi_1(g) + t^2\phi_2(g) + 
\dots,
$$
where each $\phi_i(g)$ denotes the $i$-th homogeneous component of the left-hand side,
considered as a polynomial in the variable $t$.

Each of the forms $\phi_i(g)$ is a $(i,i)$-form, and is called the
{\it $i$-th Chern form} of the metric $g$. It is a fact (see for example
\cite{Wells} for further explanations) that the cohomology class
represented by each  $\phi_i(g)$ is independent on the metric $g$, and hence it is
a topological invariant of the manifold $M$.
These cohomology classes are called the {\it Chern classes} of $M$
and they are going to be denoted by $c_i(M)$.

{\bf Remark:} We can define more generally the curvature $\Om(E)$ of a
hermitian metric $h$ on
a general complex vector bundle $E$ on a complex manifold $M$.

Let $\nabla = \nabla(h)$ be a connection on a vector bundle $E \rightarrow M$.
Then the {\em curvature} $\Om_E(\nabla)$ is defined to be the element
$\Om \in \Om^2(M, \text{End}(E, E))$ such that the $\CC$-linear mapping
$$
\Om : \Gamma(M, E) \rightarrow \Om^2(M, E)
$$
has the following representation with respect to a frame $f$:
$$
\Om(f) = \Om(\nabla, f) = d\theta(f) + \theta(f) \wedge \theta(f).
$$
Here, $\Gamma(M, E)$ is the set of sections of the vector bundle $E$,
$\Om^2(M, E)$ is the set of $E-$valued $2$-forms, and
$\theta(f)$ is the connection matrix associated with $\nabla$ and $f$
(with respect to $f$, we can write $\nabla = d + \theta(f)$).

Similarly, one defines
the Chern class $c_i(M, E)$ of a vector bundle and these are also independent on the
choice of the connection.
In fact, we use the expression ``Chern classes $c_i(M)$ of the manifold $M$''
meaning the Chern classes $c_i(M, TM)$
of the tangent bundle of M.

We will restrict our attention to the first Chern class $c_1(M)$ of the manifold $M$.
Note that the form $\phi_1(g)$
represents  the class $c_1(M)$ (by definition),
and that $\phi_1(M)$ is simply the trace of the
curvature form:
\begin{equation}
\label{e.phi1}
\phi_1(g) =  \frac{\sqrt{-1}}{2\pi} \sum_{i=1}^n \Om_i^i
= \frac{\sqrt{-1}}{2\pi} \sum_{i, p = 1}^n g^{i\pbar}  R_{i\pbar k \lbar} dz^k \wedge d\zbar^l.
\end{equation}

On the other hand, notice that the right-hand side of (\ref{e.phi1}) is equal to
$\frac{\sqrt{-1}}{2\pi} \Ric_{k\lbar}$, in view of (\ref{e.Riccidef}).
Therefore, we conclude that the Ricci form of a \kahler metric represents
the first Chern class of the manifold $M$.A natural question that arises is: given a \kahler 
class
$[\om] \in H^2(M, \RR) \cap H^{1,1}(M, \CC)$ in a compact, complex manifold
$M$, and any $(1, 1)$-form $\Om$ representing $c_1(M)$, is that possible to
find a metric $g$ on $M$ such that $\Ric(g) = \Om$?

This question was addressed to by Calabi in 1960, and it was answered by Yau \cite{Y1} almost
$20$ years later.

\begin{thm}
{\bf (Yau, 1978)}
If the  manifold $M$ is compact and \kahler, then there exists a unique (up to constant)
\kahler metric $g$ on $M$ satisfying $\Ric(g) = \Om$.
\end{thm}

This theorem has a large number of applications in different
areas of Mathematics and Physics. Its proof is based on translating
the geometric statement into a non-linear
partial differential equation, as follows.

First, fix a \kahler form $\om \in [\om]$ representing the previously given
\kahler class in $H^2(M, \RR) \cap H^{1,1}(M, \CC)$.
In local coordinates, we can write
$\om$ as $\om = g_{i\jbar}dz^i \wedge d\zbar^j$.

The $(1, 1)$-form $\Om$ is a representative for  $c_1(M)$, and
we have seen that $\Ric(\om)$
represents the same cohomology class as $\Om$.
Therefore, since $\Ric(\om)$  is also a  $(1, 1)$-form,
we have that,  due to the famous $\partial \bar \partial$-Lemma,
there exists a function $f$ on $M$ such that
$$
\Ric(\om) - \Om =  \frac{\sqrt{-1}}{2\pi} \partial \bar\partial f,
$$
where $f$ is uniquely determined after imposing the normalization
\begin{equation}
\label{e.integrability}
\int_M \left( e^f - 1 \right) \om^n = 0.
\end{equation}
Notice that $f$ is fixed once we have fixed $\om$ and $\Om$.

Again by the $\partial \bar \partial$-Lemma,
any other $(1, 1)-$form in the same cohomology class $[\om]$
will be written as
$\om + \frac{\sqrt{-1}}{2\pi} \partial \bar\partial \phi$, for some
function $\phi \in C^\infty(M, \RR)$.

Therefore, our goal is to find a representative
$\om + \frac{\sqrt{-1}}{2\pi} \partial \bar\partial \phi$ of the
class $[\om]$ that satisfies
\begin{equation}
\label{e.ric}
\Ric\left( \om + \frac{\sqrt{-1}}{2\pi} \partial \bar\partial \phi \right)
=
\Om =
\Ric(\om) - \frac{\sqrt{-1}}{2\pi} \partial \bar\partial f.
\end{equation}

Rewriting (\ref{e.ric}) in local coordinates, we have
$$
- \partial \bar\partial \log \det \left( g_{i\jbar} +
\frac{\partial^2 \phi}{\partial z_i \partial \zbar_j}
                                                  \right)
=
- \partial \bar\partial \log \det \left( g_{i\jbar}\right) - \partial \bar\partial f,
$$
or
\begin{equation}
\label{e.welldef}
 \partial \bar\partial \log
\frac{\det \left( g_{i\jbar} +
\frac{\partial^2 \phi}{\partial z_i \partial \zbar_j}
                                                  \right)}{\det \left( g_{i\jbar}\right)}
= \partial \bar\partial f.
\end{equation}

Notice that, despite of the fact that this is an expression given in
local coordinates, the term at the right-hand side of (\ref{e.welldef})
is defined globally.
Therefore, we obtain an equation
well-defined on all of $M$. In turn, this equation gives rise to the
following (global) equation
\begin{equation}
\label{e.MA100}
\left(\om + \frac{\sqrt{-1}}{2\pi}\partial \bar\partial \phi \right)^n = e^f \om^n.
\end{equation}

We shall also require positivity of the resulting \kahler form:
$\left(\om + \partial \bar\partial \phi \right) > 0$ on $M$.
This equation is a non-linear partial differential equation of
Monge-Amp\`ere type, that is going to be referred to from now on as the Complex
Monge-Amp\`ere Equation.

We remark that, if $\phi$ is a solution
to (\ref{e.MA100}),  $\om + \partial \bar\partial \phi$
is the \kahler form of our target metric $g$, {\em ie},
$\Ric(g) = \Om$.
Therefore,  in order
to find metrics that are solutions to Calabi's problem, it suffices to determine a solution
$\phi$ to (\ref{e.MA100}).

The celebrated Yau's Theorem in \cite{Y1} determines a unique solution to (\ref{e.MA100})
when $f$ satisfies the integrability condition
(\ref{e.integrability}). The proof of this result is based on the continuity method,
and we sketch here a brief outline of the proof.

The uniqueness part of Calabi Conjecture was proved by Calabi in the $50$'s.
Let $\om', \om'' \in [\om]$ be representatives of the \kahler class $[\om]$
such that $\Ric(\om') = \Ric(\om'') = \Om$. Without loss of generality, we may assume
that $\om'' = \om$, and hence $\om' = \om + \partial \bar \partial u$.

Notice that
\begin{equation}
\label{e.blah}
0 = \frac{1}{\Vol_\om(M)} \int_M u((\om')^n - \om^n) =\frac{1}{\Vol_\om(M)}
\int_M -u \partial \bar \partial u \wedge
\left[  (\om')^{n-1} + (\om')^{n-2}\wedge \om + \cdots + \om^{n-1}   \right].
\end{equation}

\noindent However $\om'$ is a \kahler form, so that $\om'>0$. We then conclude that
the right-hand side of (\ref{e.blah}) is bounded from below by
$\frac{1}{\Vol_\om(M)} \int_M -u \partial \bar \partial u \wedge w^{n-1}$.
Therefore,
\begin{equation}
0 \geq \frac{1}{\Vol_\om(M)} \int_M -u \partial \bar \partial u \wedge w^{n-1}
= \frac{1}{n\Vol_\om(M)} \int_M |\partial u|^2 w^{n}
= \frac{1}{2n\Vol_\om(M)} \int_M |\nabla u|^2 w^{n},
\end{equation}
implying that $|\nabla u| = 0$, hence $u$ is constant, proving the uniqueness
of solution to (\ref{e.MA100}).

\bigskip

Let us now consider the existence of solution to  (\ref{e.MA100}).
Define, for all $s \in [0,1]$, $f_s = sf + cs$, where
the constant $c_s$ is defined by the requirement that $f_s$ satisfies the
integrability condition $\int_M[e^{f_s} - 1 ]\om^n = 0 $.

Consider the family of equations
\begin{equation}
\label{e.fs}
\left( \om + \partial \bar \partial u_s\right)^n = e^{f_s} \om^n.
\end{equation}
We already prove that the solution $u_s$ to (\ref{e.fs}) is unique, if
it exists.

Let $A = \{s \in [0,1];$ (\ref{e.fs}) is solvable for all $t \leq s \}$.
Since $A \neq \emptyset$ ($0 \in A$), we just need to show that $A$ is open and closed.

{\bf Openness:} Let $s \in A$, and let $t$ be close to $s$. We want to show that
$t \in A$. In order to do so, let $\om_s = \om + \partial \bar \partial u_s$, for $u_s$ a
solution to  (\ref{e.fs}).
We define the operator $\Psi = \Psi_s$ by
$$
\Psi: X \rightarrow Y; \hspace{2cm}
\Psi(g) = \log\left(\frac{(\om_s + \partial \bar \partial g)^n}{w_s^n}\right),
$$
where $X$ and $Y$ are subsets (not subspaces) of
$C^{2, 1/2}(M.\RR)$ and $C^{0, 1/2}(M.\RR)$
satisfying some extra non-linear conditions.

The linearization of $\Psi$ about $g = 0$ is simply the
metric laplacian with respect to the metric $\om_s$.
By the Implicit Function Theorem, the invertibility of the
laplacian (a classical result that can be found in \cite{GT}, for example)
establishes the claim.

{\bf Closedness:}
The proof that $A$ is closed is a deep result, involving complicated
{\em a priori} estimates. A reference for this proof is Yau's paper itself \cite{Y1},
or for a more detailed proof, the books \cite{Tianbook} and \cite{Asterisque}.

\medskip

 Yau's Theorem provided a satisfatory answer to the problem of
finding Ricci-flat metrics when  the underlying manifold $M$ is compact.

Calabi's Problem, though, has a very natural generalization for the case
of a special class of open manifolds. This is the subject of the next section.

\section{Noncompact version of the Calabi-Yau Problem - Existence}
\label{s.TY}

To discuss the Calabi Conjecture on open manifolds, some minor
modifications in the original problem are needed.
Namely, we need to understand which line bundle associated to an open
manifold should play the role of the canonical line bundle $K_M$ for a compact
manifold on the Calabi-Yau problem.

Suppose that $\Mbar$ is a compact, \kahler manifold, and let $D$ be a smooth
divisor in $\Mbar$.
We are now interested in constructing complete \kahler metrics with prescribed Ricci curvature
on the open manifold $M$, defined as the complement of the divisor $D$ in $\Mbar$.

If $g'$ is a metric defined on $\Mbar$, then the metric $\det(g')$ (given locally by
$\det (g') = \det(g'_{i\jbar}) dz_1 d\zbar_1 \dots dz_n d\zbar_n$) is a metric defined on
the canonical line bundle $K_{\Mbar}$ of $\Mbar$.

Consider the line bundle $L_D$ associated to $D$, let $S$ be a defining section
of $D$ in $L_D$, and finally, let $h$ define a hermitian metric on $L_D$. Let us write
$h$, in the previous choice of local coordinates, as a positive function $a$.

With the preceding notations, the
line bundle given by $K_{\Mbar} \otimes - L_D$ has a metric
defined locally by $\det (g'_{i\jbar}) a^{-1}$. Indeed the reader will note
that this expression makes sense globally on $M$.
In turn, the metric $\det (g'_{i\jbar}) a^{-1}$ can be written
as $\det (g'_{i\jbar}) a^{-1} = \det (\frac{g'_{i\jbar}}{b})$, where $b^n = a$. In
particular, we have a new metric $g$ defined on $M$ (and also on $\Mbar$) which
is given in local coordinates by $g_{i\jbar} = \frac{g'_{i\jbar}}{b}$. Naturally
the Ricci form of the metric $g$ is a representative of the first Chern class
$c_1(-K_{\Mbar} \otimes - L_D)$. On the other hand, we would like the resulting
metric $g$ to be complete on the open manifold $M$. Strictly speaking, this will
never happen since $g$ is also a metric on the closure $\Mbar$. Nonetheless, this
construction suggests a natural way to try to obtain complete metrics. Namely
we let the metric $h$ conveniently degenerate on the divisor $D$. This implies that the
function $a$ will vanish on $D$ and thus that the metric $g$ will become unbounded
near $D$. So we may hope to find complete metrics on $M$ by this procedure. Note also
that the class of the Ricci form of $g$ is not affected by the ``degeneration'' of $h$.

Summarizing what precedes, to generalize Calabi's Conjecture to open manifolds, we begin
by fixing a representative $\Om \in c_1(-K_{\Mbar} \otimes - L_D)$.
From our previous discussion, the Ricci form of
a \kahler metric defined on $M$ is a representative of $c_1(-K_{\Mbar} \otimes - L_D)$.
Now we want to study the converse problem, namely:

\noindent {\bf Question:} Fixed a \kahler class $[\om]$ in the manifold $M$,
pick any representative $\Om$ of the first Chern class $c_1(-K_{\Mbar} \otimes - L_D)$.
Can we construct a complete \kahler metric $g$ on M such that
$\Ric(g) = \Om$?

In order to state the existence result by Tian and Yau, we need the following definitions,
that can also be found in \cite{TY1}.

\begin{defin}
Let $(X, g)$ be a complete, Riemannian manifold, and
let $K, \alpha, \beta$ be nonnegative numbers. The manifold
$(X, g)$ is of {\em $(K, \alpha, \beta)$-polynomial growth} if
\begin{itemize}
\item
The sectional curvature of $(X,g)$ is bounded by $K$;
\item
For some fixed point $x_0 \in X,$ there is a constant $C$ such that
${\text Vol}_g(B_R(x_0)) \leq C R^\alpha$ for all $R>0$,and
\item
${\text Vol}_g(B_1(x)) \geq C^{-1} (1 + d_g(x_0, x))^{-\beta}$,
\end{itemize}
where $B_R(x)$ is  the geodesic ball of radius $R$ around $x$ and
${\text Vol}_g, d_g$ denote, respectively,  the volume and the distance function
 associated to the metric $g$.
\end{defin}

\medskip

\begin{defin}
Let $(M, g)$ be a complete \kahler manifold. We say that
$(M, g)$ is of {\em quasifinite geometry of order $\ell + \de$}
if there exist $r>0$, $r_1>r_2 >0$ such that for all $x \in M$,
there exists a holomorphic map
$$\phi_x : U_x \subset (\CC^n, 0) \rightarrow B_r(x)$$ such that
\begin{itemize}
\item
$\phi_x(0) = x$, and $B_{r_2} \subset U_x \subset B_{r_1}$, where
$B_t = \{z \in \CC^n; |z| < t \}$;
\item
The pullback metric $\phi_x^\star g$ is a \kahler metric
on $U_x$ such that the metric tensor of
$\phi_x^\star g$ and its derivatives up to order $\ell$ are bounded
and $\de$-H\"older-continuously bounded with respect with the natural
coordinate system on $\CC^n$.
\end{itemize}
\end{defin}

We note that if $(M, g)$ is a \kahler manifold such that
its sectional curvature and covariant derivative of the scalar
curvature are bounded, then $(M, g)$ is of quasi-finite geometry of
order $2 + 1/2$ (as in \cite{TY1}).

\begin{thm}
\label{t.TianYau1}
{\bf (Tian, Yau, 1990)}
Let $(M, g)$ be a complete \kahler manifold of quasi-finite geometry of
order $2 + 1/2$ and with $(K, 2, \beta)$-polynomial growth.
Let $f$ be a smooth function satisfying the integrability condition
$\int_M (e^f -1)\om_g = 0$, and for some constant $C$,
\begin{itemize}
\item
$\sup\{|\nabla_g f|, |\Delta_g f|\} \leq C$;
\item
$|f(x)| \leq \frac{C}{(1 + \rho(x))^\beta}$,
\end{itemize}
where $\rho(x) = d_g(x_0, x)$ for a fixed point  $x \in M$.Then there exists a bounded, smooth 
solution $u$ to the Monge-Amp\`ere
equation
\begin{eqnarray}
\label{e.MA}
\left(\om + \frac{\sqrt{-1}}{2\pi}\partial \bar\partial u \right)^n &=& e^f \om^n, \\
\nonumber
\om + \frac{\sqrt{-1}}{2\pi}\partial \bar\partial u &>& 0 \hspace{1cm} \text{on $M$}.
\end{eqnarray}
\end{thm}

{\bf Sketch of Proof:}
The idea of the proof of this theorem is based
on a perturbation method. The function $f$ is
replaced by a sequence $\{f_m\}$ of compactly supported
functions that converge uniformly to $f$, and we consider
the sequence of modified Monge-Amp\`ere Equations
\begin{eqnarray}
\label{e.MAm}
\left(\om + \frac{\sqrt{-1}}{2\pi}\partial \bar\partial u_m \right)^n &=& e^{f_m} \om^n, \\
\nonumber
\om + \frac{\sqrt{-1}}{2\pi}\partial \bar\partial u_m &>& 0 \hspace{1cm} \text{on $M$}.
\end{eqnarray}

The strategy of solving (\ref{e.MA}) is to show that, for each $m$,
(\ref{e.MAm}) admits a solution $u_m$, and that a subsequence of $\{u_m\}$
converges uniformly to a solution $u$ of (\ref{e.MA}).

And, in order to show solvability of (\ref{e.MAm}), Tian and Yau describe
the solution $u_m$ as the uniform limit of solutions $u_{m,\ve}$ to
\begin{eqnarray}
\label{e.MAmve}
\left(\om + \frac{\sqrt{-1}}{2\pi}\partial \bar\partial u_{m,\ve} \right)^n &=& e^{f_m + \ve  
u_{m,\ve} } \om^n, \\
\nonumber
\om + \frac{\sqrt{-1}}{2\pi}\partial \bar\partial  u_{m,\ve}  &>& 0 \hspace{1cm} \text{on $M$}.
\end{eqnarray}

The solvability of (\ref{e.MAmve}) is due to Cheng and Yau \cite{CY1}. By standart elliptic
theory, the main step in showing that $u_{m, \ve} \rightarrow u_m$ is the uniform $C^0$ 
estimate for
the solutions $u_{m, \ve}$ (some extra work will give the $C^{2, 1/2}$ uniform estimates, which 
are
sufficient to guarantee the result). In order to derive such estimates, a key step is the 
derivation
of weighted Sobolev inequalities. Further detail in the proof of this theorem can
be found in the paper \cite{TY1}.
\qed

This theorem, though, does not provide enough information
about the behavior of the solution $u$ close to the divisor.
This is an interesting question to be studied, with a number of applications.
In order to illustrate the applications, let us quickly discuss how do they
enter into the construction of an example due to Kovalev \cite{Ko} of a
Riemannian manifold having special holonomy group $G_2$. Here we remind the
reader that examples of Riemannian manifolds with holonomy $G_2$ are hard
to construct though these manifolds play an important role
in Physics.

\bigskip

\noindent {\bf  Application: Manifolds with Special Holonomy $G_2$ (Kovalev, \cite{Ko})}.

The  manifold in question has
dimension~$7$ and will be obtained by means of a connected-sum construction.

The strategy is to glue appropriately two copies of
$S^1 \times W$, where $W$ is a $6$-manifold  with
Holonomy group $S\cal U$$(3)$, a maximal subgroup of $G_2$.

An example of such $W$ would be a complete \kahler manifold with
zero Ricci curvature and asymptotically cylindrical end, diffeomorphic
to $D\times \RR_+ \times S^1$, where $D$ is a simply-connected, compact
\kahler manifold with $c_1(D) = 0$.

Hence, a main step in Kovalev's construction is to
find a Ricci-flat metric on $W$ that is asymptotically
close to a {\em product} Ricci-flat metric on
$D\times \RR_+ \times S^1$.
Once this is done, it is possible to show (\cite{Ko})
that this metric has holonomy group $S\cal U$$(3)$.

The Riemannian product of
$W$ and $S^1$ yields open $7$-manifolds with same holonomy
group as $W$. Finally, in order to obtain the manifold with
holonomy group $G_2$, it is necessary to glue (in a non-standart way)
two copies of $S^1 \times W$, interchanging the $S^1$-factors
(recall that the end of $W$ fibers over $D \times \RR_+$) in order
to avoid infinite fundamental group, a necessary condition for
a manifold to admit a metric with
holonomy group $G_2$.

The existence of a Ricci-flat metric in this context  is attributed to
Tian and Yau (\cite{TY1}), but the case dealt with in the
paper \cite{TY1} is different, and doesn't apply directly to the case above.
The asymptotics of this problem can be (partially) derived from
the results in Section~\ref{s.general} of this paper, and in a future
work, we will present a detailed proof of the existence result in the
context of complete \kahler manifold with
zero Ricci curvature and asymptotically cylindrical end, as well as its
detailed asymptotics.

\section{Construction of complete \kahler metrics with special approximating properties}

The goal of this section is to discuss recent developments
on the study of the behavior of complete Ricci-flat metrics at infinity.

In the authors' thesis (\cite{Santoro1}, \cite{Santoro2}), the strategy to
study the asymptotics of Tian-Yau solution to (\ref{e.MA}) is divided in two
main steps.

The first of the is the inductive construction of an explicitly given sequence of
complete \kahler metrics on the manifold $M$ with special approximating properties.The other 
step is to use these metrics as ambient metrics, and study the solution
to (\ref{e.MA}) on $M$.

To state the main results of \cite{Santoro2}, let us consider a
compact, complex manifold $\Mbar$ of complex dimension $n$. Let
$D$ be an {\em admissible} divisor in  $\Mbar$, ie, a divisor satisfying the
following conditions:
\begin{itemize}
\item
Sing{ $\Mbar$}$\subset D$.
\item
$D$ is smooth in $\Mbar\setminus $ Sing{ $\Mbar$}.
\item
For every $x \in$ Sing{ $\Mbar$}, the corresponding
local uniformization $\Pi_x : \tilde{\cal U}_x \rightarrow
{\cal U}_x$, with $\tilde{\cal U}_x \in \CC^n$, is such that
$\Pi_x^{-1}(D)$ is smooth in $\tilde{\cal U}_x $.
\end{itemize}

Let $\Om$ be a  smooth, closed $(1,1)$-form in the
cohomology class $c_1(K_{\Mbar}^{-1} \otimes L_D^{-1})$.
Let $S$ be a defining section of $D$ on $L_D$
and let $M$ be the open manifold $M = \Mbar \setminus D$.
Consider a hermitian metric $||.||$ on $L_D$.

\begin{thm}
\label{t.approximatingmetrics}
{\bf \cite{Santoro2}}
Let $M$, $\Om$ and $D$ be as above. Then, for every $\ve >0$, there exists an
explicitly given complete
\kahler metric $g_\ve$ such that
\begin{equation}
\label{e.MAve}
\Ric(g_\ve) - \Om = \frac{\sqrt{-1}}{2\pi} \partial \bar\partial f_\ve
\hspace{1.5cm} \text{on $M$,}
\end{equation}
where $f_\ve$ is a smooth function on $M$ that decays to the order of $O(||S||^\ve)$.
Furthermore, the Riemann curvature tensor $R(g_\ve)$ of the metric $g_\ve$ decays
to the order of $O((-n \log ||S||^2)^{-\frac{1}{n}})$.
\end{thm}

Note that the metrics given above are explicitly given, and
its construction is an interesting result in itself.

To prove  Theorem \ref{t.approximatingmetrics}, we start by fixing
an orbifold hermitian metric $||.||$ on $L_D$ such that its curvature form
$\tilde\om$ is positive definite along $D$. Here we are using the assumption that
$D$ is ample. In Section~\ref{s.general}, we discuss how we can weaken this assumption
on the divisor, providing a generalization on the statements in \cite{Santoro2}.

Yau's Theorem for compact manifolds implies that,
by rescaling $||.||$ by an appropriate factor, we may assume
that $\tilde\om$, when restricted to the
infinity $D$, defines a metric $g_D$ such that $\Ric(g_D) = \Om|_D$.
Denote by $||.||_\phi  = e^{-\phi/2} ||.||$  the rescaling
of the metric $||.||$, where $\phi$ is any smooth function on $\Mbar$.

We define
\begin{equation}
\label{omegaphi}
\om_\phi = \frac{\sqrt{-1}}{2\pi} \frac{n^{1 + 1/n}}{n+1} \partial \bar{\partial}
(-\log||S||_\phi^2)^{ \frac{n+1}{n}}.
\end{equation}

The key step on the proof of Theorem \ref{t.approximatingmetrics} is its local version,
as follows.
\begin{prop}
\label{t.inductive}
Let $\Mbar$ be a compact \kahler manifold of complex dimension $n$, and let $D$ be
an admissible divisor in  $\Mbar$. Consider also
a form  $\Om \in c_1(-K_{\Mbar} - L_D)$, where $L_D$ is
the line bundle induced by $D$.

Then there exist sequences of neighborhoods $\{V_m\}_{m \in \NN}$ of $D$
along with complete K\"ahler
metrics $\om_m$  on $(V_m\setminus D, \partial(V_m\setminus D))$ (as defined in \ref{omegaphi})
such that
\begin{equation}
\label{e.thmseqs}
\text{Ric}(\om_m) - \Om = \frac{\sqrt{-1}}{2\pi}\partial \bar{\partial}f_m  \\ \text{
on $V_m\setminus D$}
\end{equation}
where $f_m$ are smooth functions on $M = \Mbar \setminus D$. Furthermore
each $f_m$ decays to the order of $O(||S||^m)$. In addition
the curvature tensors $R(g_m)$ of the metrics  $g_m$
decay at least to the order of $(-n\log||S||_m^2)^{\frac{-1}{n}}$ near the divisor.
\end{prop}

Some further work (\cite{Santoro2}) allow us to extend each of these metrics accordingly to the
whole manifold.

Let us say a few words about the proof of Proposition~\ref{t.inductive}.
Fefferman, in his paper \cite{Fef}, developed inductively an
\hbox {$n$-th} order approximation to a complete
\KE metric on strictly pseudoconvex domains on
$\mathbb{C}^n$ with smooth boundary, and he suggested that higher order
approximations could be
obtained by considering $\log$ terms in the formal expansion of
the solution to a certain complex Monge-Amp\`ere
equation. This idea was taken up by Lee and
Melrose in \cite{Mel}, where they constructed the full
asymptotic expansion of the solution to the
Monge-Amp\`ere equation considered by Fefferman, determining completely
the form of the singularity and improving  the regularity of the existence result of
Cheng and Yau \cite{CY1}.

Tian and Yau \cite{TY3} showed the existence of \KE metrics under certain conditions
on the divisor, and Wu \cite{DaminWu} developed the asymptotic expansion to the
Cheng-Yau metric on a quasi-projective manifold (also assuming some conditions
on the divisor), as the parallel part to the work
of Lee and Melrose \cite{Mel}.

As mentioned in the Introduction, the K\"ahler-Einstein condition used in these
works makes the underlying analysis very special and widely developed by many people,
and does not apply to the Ricci-flat case.
Yet the proof of  Proposition~\ref{t.inductive} is inspired on
the inductive methods
of Fefferman and Lee-Melrose, applied to  the context of complete Ricci-flat metrics
quasi-projective manifolds. The full detailed version of it can be found in
\cite{Santoro1}, \cite{Santoro2}.

\noindent
{\bf Idea of the proof of  Proposition~\ref{t.inductive}:}

First, consider the function
$$
f_\phi(x)= -\log||S||^2  - \log(\frac{\om_\phi^n}{{\om'}^n}) - \Psi,
$$
where $\om'$ is any \kahler form in $\Mbar$, and $\Psi$ is related to $\om'$
by $\Om = \text{Ric}(g') - \tilde\om +  \frac{\sqrt{-1}}{2\pi} \partial \bar{\partial} \Psi.$

When $\phi = 0$, it
is not hard to see that
the function $f_0(x)$ will converge uniformly to a constant if and only if
$\Ric(g_D) = \Om|_D$. Without loss, we can assume that this constant is zero.

It is possible to extend $f_0$ smoothly to  be zero along the divisor.
Hence, there exists a $\delta_0 > 0$ such that, in the neighborhood
$V_0 := \{x \in \Mbar ; ||S(x)|| <  \delta_0\}$,
$f_0$ can be written as
$$
f_0 = S \cdot u_1 + \Sbar \cdot \bar{u}_1,
$$
where $u_1$ is a $C^\infty$ local section in $\Gamma(V_0, L_D^{-1})$.

Our goal now would be to construct a function $\phi_1$ of the form
$S \cdot \theta_1 + \Sbar \cdot \overline{\theta}_1$, so that the
corresponding $f_{\phi_1} = f_1$ vanishes at order $2$ along $D$, and then
proceed inductively to higher order. Unfortunately, there is an obstruction
to higher order approximation that lies in the kernel of the laplacian on
$L_D^{-1}$ restricted to $D$.
In order to deal with this difficulty, we must introduce $(-\log||S||^2)$-terms in
the expansion of $\phi_1$.

Following the techniques in \cite{TY2}, we construct
inductively in \cite{Santoro2} a sequence of hermitian metrics $\{||.||_m\}_{m>0}$ on
$L_D$ such that, for any $m>0$, there exists a $\de_m >0$ satisfying:
\begin{enumerate}
\item
The corresponding K\"ahler form $\om_m$ associated to $||.||_m$ (as defined in
(\ref{omegaphi})) is positive definite in $V_m := \{x \in \Mbar ; ||S(x)|| <  \delta_m\}$.
\item
The function $f_m$ associated to $\om_m$ can be expanded
in $V_m$ as
\begin{equation}
\label{expansionfm}
f_m =
\sum_{k\geq m+1} \sum_{\ell =0}^{\ell_k} u_{k\ell}(-\log||S||^2_m)^\ell,
\end{equation}
where $u_{k\ell}$ are smooth functions on $\overline{V}_m$ that vanish to order
$k$ on $D$. In particular the function $u_{k\ell}$ can
be written as
$$
u_{k\ell} = \sum_{i+j = k}S^i \Sbar^j \theta_{ij} + S^j \Sbar^i \overline{\theta}_{ij},
$$
for  $\theta_{ij} \in \Gamma(V_m, L_D^{-i} \bigotimes  \overline{L}_D^{-j})$.
\end{enumerate}
This construction is very technical, and it is based on studying the kernel of a second-order
differential operator. A detailed exposition of this result can be found in \cite{Santoro2}.

\section{Asymptotics of the Monge-Amp\`ere Equation on quasi-projective manifolds}
\label{s.laplacian}

In this section, we want to conclude our second step: considering the
metrics
given by Theorem~\ref{t.approximatingmetrics} as ambient metrics,
to study the asymptotics of the solution
to the Complex Monge-Amp\`ere Equation.

\begin{thm}
\label{t.mainthm}
{\bf \cite{Santoro2}}
For each $\ve>0$, let $g_\ve$ and $f_\ve$ be given by Theorem
\ref{t.approximatingmetrics}.
As before, let $S$ be a defining section for the divisor $D$.

Consider the
solution $u_\ve$ to the problem
\begin{equation}
\label{MA1ve}
\begin{cases}
\left(
\om_{g_\ve} + \frac{\sqrt{-1}}{2\pi}\partial\bar\partial u_\ve
\right)^n  = e^{f_\ve} \om_{g_\ve}^n, & \\
\om_{g_\ve} + \frac{\sqrt{-1}}{2\pi}\partial\bar\partial u_\ve >0,
& \text{$u_\ve \in C^\infty(M, \RR)$}.
\end{cases}
\end{equation}

Then the solution $u_\ve$ decays as $O(||S||^\ve)$ near the divisor.
\end{thm}

Theorem~\ref{t.mainthm} shows that the metrics $\om_{g_\ve}$ given by
Theorem \ref{t.approximatingmetrics} are in fact a very good description
of what happens with
the actual solution $\om_{g_\ve} + \frac{\sqrt{-1}}{2\pi}\partial\bar\partial u_\ve$
at infinity.

We need to observe that, in fact, the manifold constructed in
Theorem~\ref{t.approximatingmetrics} actually satisfy the conditions on
Tian-Yau's Theorem (see \cite{Santoro1}.

Clearly, it suffices to prove the asymptotic assertions on $u_\ve$ for a small tubular
neighborhood of $D$ in $\Mbar$. Therefore, we will consider the equation
\begin{equation}
\label{e.MAomm}
\left(
\om_m +  \frac{\sqrt{-1}}{2 \pi} \partial \bar \partial u_m
\right)^n
=
e^{f_m}\om_m^n  \hspace{1cm} \text{on $V_m \setminus D$,}
\end{equation}
where $f_m$ $||.||_m$ and $\om_m$ are given by Proposition~\ref{t.inductive} and
in the discussion that follows this proposition.
We remind the reader that $|f_m|_{g_m}$ is of order of $O(||S||^m_m)$.

The proof of Theorem~\ref{t.mainthm} is based on the very careful construction
of a barrier function, that we state here for completeness.
\begin{lemma}
\label{l.barrier}
On the neighborhood \hbox{$V_m\setminus D = \{ 0< ||S||_m < \de_m \}$,}
we have
\begin{multline}
   \left\{
\om_m + \frac{\sqrt{-1}}{2\pi}\partial \bar \partial\left(C
\left[S^i \Sbar^j \theta_{ij} + \Sbar^i S^j \bar\theta_{ij} \right]
(-n\log(||S||_m^2))^k \right)
  \right\}^n
 = \\ =
\om_m^n
\left[
1  + C (-n\log(||S||_m^2))^{k-\frac{n+1}{n}}
\left\{
ij(-n\log(||S||_m^2))^2  \left[S^i \Sbar^j \theta_{ij} + \Sbar^i S^j \bar\theta_{ij} \right]
- \right. \right.\\ \left.\left. -
(-n\log(||S||_m^2)) \left[
\left( k(i+j) +j(n-1)\right) S^i \Sbar^j \theta_{ij}
+
\left( k(i+j) +i(n-1)\right)  \Sbar^i S^j \bar\theta_{ij}
\right]
\right. \right. \\ \left. \left.
+ k(k-n)
\right\}
+ O(||S||_m^{i+j+1})
\right],
\end{multline}
where  $\theta_{ij}$ is a $C^\infty$ local section of $L_D^{-i} \otimes  \overline{L}_D^{-j}$
on $V_m$.
\end{lemma}

Using Lemma~\ref{l.barrier} and the maximum principle for the
complex Monge-Amp\`ere operator, we can show:
\begin{prop}
\label{p.uvanishesalot}
Let $u_m$ be a solution to the Monge-Amp\`ere equation (\ref{MA1ve}).
If $u_m(x)$ converges uniformly to zero as $x$ approaches the divisor,
then there exists a constant $C = C(m)$ such that
\begin{equation}
\label{e.uvanishesalot}
|u_m(x)| \leq C ||S||_m^{m+1} \hspace{1cm} \text{on $V_m \setminus D$}.
\end{equation}
\end{prop}

In order to complete the proof of Theorem~\ref{t.mainthm}, we need to show that,
for a fixed $m>2$,
the solution $u_m(x)$ converges uniformly
to zero as $x$ approaches the divisor $D$.

The idea of this proof is to use again the maximum principle, with the
new extra ingredient: as pointed out before,  in \cite{TY1},
the solution $u_m$ to the Monge-Amp\`ere equation~(\ref{e.MA}) is
obtained as the uniform limit, as $\ve$ goes to {\it zero}, of solutions $u_{m,\ve}$
of the perturbed Monge-Amp\`ere equations
\begin{equation}
\label{MAve}
\begin{cases}
\left(
\om_m + \frac{\sqrt{-1}}{2\pi}\partial\bar\partial u
\right)^n  = e^{f_m + \ve u} \om_m^n, & \\
\om_m + \frac{\sqrt{-1}}{2\pi}\partial\bar\partial u >0,
& \text{$u \in C^\infty(M, \RR)$}.
\end{cases}
\end{equation}

It turns out (\cite{Santoro2}, Proposition $5.2$) that it is possible
to construct a barrier function that is {\em uniform} on $\ve$, so that
the estimates will carry out to the limiting function $u_m$.

Finally, we remark that, by the work of Cheng and Yau  \cite{CY1},
$u_{m,\ve}$ vanishes uniformly as we approach  $D$. Therefore we can apply the
maximum principle to conclude that $u_m(x)$ converges uniformly to zero as
$x$ converges to $D$.

This concludes the proof of Theorem~\ref{t.mainthm}.

\section{Further Developments}
\label{s.general}

So far, we have been considering the problem of constructing
Ricci-flat metrics on quasi-projective manifolds $M$ that can
be compactified by adding a smooth and {\em ample} divisor $D$.

The goal of this section is to show that the construction in
Theorem~\ref{t.approximatingmetrics} can be generalized to the
case where the divisor
$D$ is a {\em semi-ample} divisor. The definition of a
semi-ample line bundle/divisor is given a few paragraphs below.

Even though we still need to assume a few further conditions
on $D$,
this class of manifolds is quite large,
including in particular the asymptotically cylindrical manifolds  considered by
Kovalev (\cite{Ko}), described at the end of Section~\ref{s.TY}.

In a subsequent work, we expect to extend the existence results by Tian and Yau
(Theorem~\ref{t.TianYau1}) to asymptotically cylindrical manifolds, as well as
to the case where $D$ is allowed to have normal crossings. In those cases, we also
expect to be able to determine the asymptotic behavior of the solutions.

\bigskip

Consider a compact, complex manifold $\Mbar$ of dimension $n$, and
let $D$ be a smooth divisor on $\Mbar$.

Let $L_D \rightarrow \Mbar$ be the line bundle  associated to $D$, and for
$m \in \NN$ sufficiently large, consider the linear system
$$
|mL_D| = H^0(\Mbar, m L_D),
$$
the set of sections of $mL_D$.
Assume that $L_D$ is {\em semi-ample}, {\em ie},
the ({\em a priori} biholomorphic) map  given by
the linear system
\begin{eqnarray*}
|m L_D| &:& \Mbar \rightarrow \Mbar_{can} \subset \CC {\mathbb P}^{N_m} \\
& & p \mapsto [S_0^{(m)}(p), \cdots , S_{N_m}^{(m)}(p)]
\end{eqnarray*}
is actually a holomorphic map for $m > > 1$. Here, $\{ S_0^{(m)}(p), \cdots , S_{N_m}^{(m)}(p) 
\}$
form a basis for $H^0(\Mbar, m L_D)$.

Assuming that $L_D$ is semi-ample implies that this line bundle is also {\em semi-positive},
{\em ie}, that there exists a hermitian metric $||.||$ on $L_D$ such that
its curvature form $\tilde \om$ is semi-positive definite along $D$ ($\tilde \om \geq 0$).
Assume further that $\tilde \om$ has constant rank equals to $k$, $0 \leq k < n-1$.
This implies that the Kodaira dimension of $L_D$ is equal to $k$, {\em ie},
$\dim (\Mbar_{can}) = k$.

Consider the restriction of $|m L_D|$ to $D$, and let $D_{can} \subset \CC {\mathbb P}^{N_m}$
be its image under  $|m L_D|$. We are interested in the fibration $D \rightarrow D_{can}$, and
we will refer to the fibers of this fibrations as $D_s$, for $s \in D_{can}$.

As before, consider a $(1,1)$-form
 $\Om \in c_1(-K_{\Mbar} - L_D)$. Our goal is to construct a complete
 metric $g$ on $M = \Mbar \setminus D$ such that
$\Ric(g) - \Om = \partial \bar \partial f$,
for a function $f$ with sufficiently fast decay.

Let $\om_F$ be a $(1, 1)$- form on $D$ such that
\begin{description}
\item{1.}
$\Ric(\om_F) \in [\Om|_D]$
\item{2.}
$\om_F|_{D_s}$ is independent of $s$, for all $s$ in
$\{s \in D_{can}; D_s \text{ is non-singular} \}$.
\end{description}
We can pick $\om_F$ satisfying those two conditions, since we can use Yau's Theorem
on $D$, and $D \rightarrow D_{can}$ is a regular fibration away from the singular fibers,
so all fibers are cohomologous.

If the regular fibers of $D \rightarrow D_{can}$ were assumed to be Calabi-Yau manifolds
($c_1(D) = 0$ and simply-connected), we could define the so-called {\em semi-flat} metrics
(metrics that are Ricci-flat on regular fibers), considered by Gross and Wilson
in \cite{GW}, for the 2-dimensional case (fibers being $K3$ surfaces).

Notice that $\om_F$ is not necessarily a \kahler form anymore. In the sequel
we will define the extension of $\om_F$ (in the same cohomology
class as $\om_F$) to a neighborhood of $D \in \Mbar$. With a small abuse
of notation, this extension will still be denoted by $\om_F$.


We define, for $k = $rank$(\tilde\om)$,
\begin{equation}
\label{omegak}
\om_\phi = \frac{\sqrt{-1}}{2\pi} \frac{k^{1 + 1/k}}{k+1} \partial \bar{\partial}
(-\log||S||_\phi^2)^{ \frac{k+1}{k}}.
\end{equation}

\begin{equation}
\label{e.omegakn}
\om_\phi =
(-k \log ||S||_\phi)^{\frac{1}{k}}\tilde\om_\phi +
(-k \log ||S||_\phi)^{\frac{1-k}{k}} \frac{\sqrt{-1}}{2\pi}
\partial \log ||S||_\phi^2 \wedge \bar\partial \log ||S||_\phi^2,
\end{equation}
where $||.||_\phi = e^{\phi/2}||.||$ denotes a rescaling of the original
metric $||.||$ and $\tilde\om_\phi$ is the curvature form of the rescaled
metric $||.||_\phi$.

From (\ref{e.omegakn}), we can see that $\om_\phi$ will have rank $k$ near $D$,
since $\tilde\om_\phi$ has rank $k$ along $D$, being a simple rescaling of $\tilde \om$.

Now, we define
\begin{equation}
\eta_\phi = \om_\phi + (-k \log ||S||_\phi)^{\frac{1}{k}} \om_F,
\end{equation}
that will be locally a \kahler form due to the observations above.

Note that
\begin{equation}
\eta_\phi^n = C(k, n) \om_\phi^k \wedge (-k \log ||S||_\phi)^{\frac{n-k}{k}} \om_F^{n-k}.
\end{equation}

Define, analogously to the  ``$D$ ample'' case,
$$
f_\phi(x) = - \log(||S||)^2 - \log\left( \frac{\eta_\phi^n}{\om'}\right) - \Psi,
$$
where $\om'$ is any fixed \kahler form on $\Mbar$ and it is related to the function
$\Psi$ by $\Ric(\om') - \Om = \frac{\sqrt{-1}}{2 \pi}\partial \bar \partial \Psi$.

The key point is that the difference of two functions $f_{\phi_1}$
and $f_{\phi_2}$ will be only based on the term $\frac{\om_{\phi_1}}{\om_{\phi_2}}$,
since the part involving $\om_F$ will cancel out on the computations.
Hence, we are able to prove an analogous lemma to \cite{Santoro2}, Lemma 2.2, that is the
main technical lemma that allow us to complete
the proof of Theorem~\ref{t.approximatingmetrics}.
Hence, since the computations carry out very similarly to the
``$D$ ample'' case, we will simply state the final result, that extends one of the main
results in \cite{Santoro2}.
\begin{thm}
\label{t.extension}
Let $M$ be a quasi-projective manifold that can be compactified by adding a smooth,
semi-ample divisor $D$, that further satisfy the conditions discussed above.
 Then, for every $\ve >0$, there exists an
explicitly given complete
\kahler metric $g_\ve$ such that
\begin{equation}
\label{e.MAve2}
\Ric(g_\ve) - \Om = \frac{\sqrt{-1}}{2\pi} \partial \bar\partial f_\ve
\hspace{1.5cm} \text{on $M$,}
\end{equation}
where $f_\ve$ is a smooth function on $M$ that decays to the order of $O(||S||^\ve)$.
\end{thm}

\vskip 1cm

\flushleft

\medskip

{\bf Bianca Santoro} \ \  (bsantoro@msri.org)\\
Mathematical Sciences Research Institute\\
17 Gauss Way, room 217\\
Berkeley, CA  94720\\
USA\\

\end{document}